\documentclass{amsart}

\theoremstyle{definition}

\theoremstyle{remark}

\numberwithin{equation}{section}

\usepackage{amssymb,amsmath}

\def\s{\sigma}

\def\o{\omega}
\def\u{\mu}

\def\p{\rho}

\def\xk{x^k(\mu)}
\def\xk1{x^{k-1}(\mu)}

\def\cal#1{\mathcal{#1}}
\def\K{\cal{K}}
\def\K#1{\cal{K}_{#1}}
\def\C{\cal{C}}
\def\g{g_{\o}(\mu)}
\def\g#1{g_{#1}(\mu)}
\def\mo{{\mu}_{\o}}
\def\c{\circ}

\def\dfrac#1#2{\displaystyle{\frac{#1}{#2}}}

\usepackage{graphicx}
\graphicspath{{converted_graphics/}}

\begin{document}

\title{A Solution to the Monotonicity Problem for Unimodal Families}

\author{John Taylor}
\address{Department of Mathematical Sciences, United Arab Emirates University, Al Ain, UAE }
%\curraddr{Department of Mathematics and Statistics,
%Case Western Reserve University, Cleveland, Ohio 43403}
\email{john.taylor@uaeu.ac.ae}
\thanks{}


\date{May 15, 2008}


\begin{abstract}
In this note we consider a collection $\cal{C}$ of one parameter families of unimodal maps of $[0,1].$ Each family in the collection has the form $\{\mu f\}$ where $\mu\in [0,1].$  Denoting the kneading sequence of $\mu f$ by $K(\mu f)$, we will prove that for each member of $\cal{C}$, the map $\mu\mapsto K(\mu f)$ is monotone.  It then follows that for each member of $\cal{C}$ the map $\mu\mapsto h(\mu f)$ is monotone, where $h(\u f)$ is the topological entropy of $\mu f.$  For interest, $\mu f(x)=4\mu x(1-x)$ and $\mu f(x)=\mu\sin(\pi x)$ are shown to belong to $\cal{C}.$  
\end{abstract}

\maketitle

\section*{introduction}

%Maps on the interval provide the simplest examples to illustrate the problem of understanding how dynamic complexity evolves under deformations of a dynamical system.

Metropolis, Stein and Stein were among the first, to my knowledge, to study what are now called finite kneading sequences.  These were associated with super stable limit cycles of one parameter families of interval maps, which included $\mu f(x)=4\mu x(1-x)$ and $\mu f(x)=\mu\sin(\pi x)$ [2].  Computer studies strongly suggested a universal topological dynamics for a large class of such families, and many workers were quickly drawn to this fascinating field of study. In the 1980's and 1990's there was intense interest in the behavior of the logistic map (affinely modified) in the setting of one complex dimension.  A central question was (essentially) under what circumstances would finite kneading sequences be monotone with the parameter, as this question was associated with the structure of the Mandelbrot set, among other things.  This question was successfully addressed (in the special case of a real quadratic) in [1,5].  Here we address the question generally for a large class that includes the logistic map, and give a sufficient condition for the solution.  

All known proofs of this type apply to the case of a quadratic polynomial only and use complex analytic methods (holomorphic techniques) or depend on complex analysis (Compare [1], [5], [6], [7] [8]).  These methods are not used here.

Let $I=[0,1].$  Consider the collection of parameterized maps $\{\mu f\ |\ \mu f:I\rightarrow I\}$ with $\mu f$ at least $C^{(3)}$ in $x$.  Notice that since $\u f(x)$ is linear in $\u, $ it is $C^{\infty}$ in $\u.$  Denote the single critical point $c\in(0,1)$ and scale the map so that $f(c)=1,$ requiring that $f(0)=0=f(1).$  Then $\mu f(c)=\mu.$  Denote the $n^{th}$ iterate of $\mu f$ by $f_\mu^n(x)=(\mu f)\circ\cdots\circ(\mu f)\circ(\mu f)(x),$ where the composition is $n$-fold.
\bigskip

For any $x\in I$, the {\it orbit} of $x$ is the set $O(x)=\{f_\mu^n(x)| n\ge 0\}.$  Associate with $O(x)$ the word   
$\o(x) = \o_0\o_1\o_2\cdots$ with $\o_k\in\{L, C, R\}$ where words are formed as follows:

$$\o_k = \left\{ \aligned L, \quad \text{for} \quad f_\mu^k(x) < c\cr
C, \quad \text{for} \quad f_\mu^k(x) = c\cr
R, \quad \text{for} \quad f_\mu^k(x) > c \endaligned \right.$$
\bigskip

\noindent $\o(x)$ is called the {\em itinerary} of $x$ under $f$.
We are interested in studying the itinerary associated with $O(\u)$. This special itinerary is called the {\it kneading sequence} of $\mu f$, symbolized $K(\u f)=\o(\u)$.   
\bigskip

In particular, we concentrate on finite itineraries having the form $\o_0\o_1\o_2\cdots C,$ as these correspond to (super stable) periodic points.  
\bigskip

%%%%%%%%%%%%%%%%%%%%%%%%%%%%%%%%

Define $\cal{K}=\{\o\ |\ 1\le|\o|<\infty,\ \o=K(\mu f),\ f\in\cal{C}\}$.  That is, $\cal{K}$ is exactly all finite kneading sequences for any member of $\cal{C}.$

Let $\cal{K}_n$ represent all kneading sequences of length $n.$  For each $n\ge 1,$ it is known that
 $$|K_n|=\displaystyle{\dfrac{1}{2n}\sum\mu(d)2^{n/d}}$$
where the sum is taken over all odd square free divisors of $n.$ [5]
\bigskip
%%%%%%%%%%%%%%%%%%%%%%%%%%%%%%%%

The following preliminaries are necessary to give meaning to the statement $\mu\mapsto K(\mu f)$ is monotone;  explicitly, we need a total order $\prec$ on the kneading sequences that reflects the
order of the real line in the sense that
 $x < y$ implies $\o(x) \preccurlyeq \o(y)$.
\bigskip

%\noindent {\bf Note.}  
This is done by defining $L\prec C\prec R.$  Then, if
$A=\{a_k\}\ne B=\{b_k\},$ let $N$ be the smallest index for which $a_N\ne b_N$
and let $\p_{N-1}$ be the number of $R's$ in the word $a_1\cdots a_{N-1}.$
Then define $A\prec B$ if $a_N\prec b_N$ and $\p_{N-1}$ is even or if $a_N\succ b_N$ and $\p_{N-1}$ is odd.

This order it sometimes referred to as the {\em parity-lexicographic} order.  
The intuition derives from the fact that  $\u f$ is orientation preserving or reversing according as $x\in(0,c)$ or $x\in(c,1).$ 

A word $\o$ is called {\it maximal} (or {\it shift-maximal}) provided it
is greater (in the parity lexicographic order) than all of its shifts,
where, as usual, the shift operator $\sigma$ is defined by the
action $\sigma(\o) = \o_1\o_2 \o_3\cdots$ on the word $\o=\o_0\o_1\o_2\o_3\cdots .$  
\bigskip

In kneading theory there are several versions of an ``intermediate value
theorem''.  This type of theorem is fundamental in that it relates abstract words to the behavior of dynamical systems.  That is, it connects the set of kneading sequences ordered by the relation $\prec$ and the parameter space (an interval in the real line)
with the usual order.  The following version is essentially that found in
[3]: 
\bigskip

{\bf Theorem A}\ Let $\{\mu f \}$ be any one parameter family of
$C^1$ unimodal maps.  If $\mu_1 < \mu_2$ are two parameter values with
corresponding kneading sequences $K(\mu_1 f)\prec K(\mu_2 f)$, and if
$\o$ is any shift-maximal sequence with the property that
$$K(\mu_1 f)\prec \o\prec K(\mu_2 f),$$
then there exists a $\mu$ such that $\mu_1 < \mu < \mu_2$ and $\o =
K(\mu f)$.
\vfill\eject

\begin{center}
{\large\textsc{Iterating Inverse Functions Along a Word}}
\end{center}
\bigskip

Since $\mu f$ double covers $[0, \u]$ in such a way that $\mu f([0,c])=[0, \u]=\mu f([c,1]),$  the functions $(\mu f)_{\o}^{-1},\ \o\in\{L, R\}\ \text{ have the action}$
$$ (\mu f)_{L}^{-1}([0, \u])=[0,c]\ \text{and}\ (\mu f)_{R}^{-1}([0, \u])=[c,1].$$
\smallskip

Notice that, given $f$, we can write $(\mu f)_{\o}^{-1}(x)$ as a function of three variables, one of which is a word:
$$
G_{\o}(\mu, x)=(\mu f)_{\o}^{-1}(x)
$$

Let $\o=\o_{1}\o_{2}\cdots\o_{n}.$
\bigskip

Define $G_{\o_n}(\mu, x)=(\mu f)^{-1}_{\o_n}(x)|_{x=c}=(\mu f)^{-1}_{\o_n}(c)=g_{\o_n}(\mu);$
\bigskip

$G_{\o_{n-1}\o_n}[\mu, x]=G_{\o_{n-1}}[\mu, g_{\o_n}(\mu)]=G_{\o_{n-1}}[\mu, (\mu f)^{-1}_{\o_n}(c)]$
\bigskip

$=(\mu f)^{-1}_{\o_{n-1}}\left[(\mu f)^{-1}_{\o_n}(c)\right]$
\bigskip

$=(\mu f)^{-1}_{\o_{n-1}}\circ(\mu f)^{-1}_{\o_{n}}(c)$
\bigskip

$=g_{\o_{n-1}}\circ g_{\o_{n}}(\mu)$
\bigskip

$=g_{\o_{n-1}\o_{n}}(\mu).$
\bigskip

Continuing,
\bigskip

$G_{\o_{1}\o_{2}\cdots\o_{n}}\left[\mu, x\right]=G_{\o_{1}}\left[\mu, g_{\o_{2}\cdots\o_{n}}(\mu)\right]$
\bigskip

$=G_{\o_{1}}[\mu, (\mu f)^{-1}_{\o_2}\circ \cdots \circ (\mu f)^{-1}_{\o_{n}}(c)]$
\bigskip

$=(\mu f)^{-1}_{\o_1}\left[(\mu f)^{-1}_{\o_2}\circ \cdots \circ (\mu f)^{-1}_{\o_{n}}(c)\right]$
\bigskip

$=(\mu f)^{-1}_{\o_1}\circ (\mu f)^{-1}_{\o_2}\circ \cdots \circ (\mu f)^{-1}_{\o_{n}}(c)$
\bigskip

$=g_{\o_{1}}\circ g_{\o_{2}}\circ \cdots \circ g_{\o_{n}}(\mu)$
\bigskip

$=g_{\o_{1}\o_{2}\cdots \o_{n}}(\mu) $
\bigskip

$=g_{\o}(\mu). $
\bigskip

Observe that the iteration is in the variable  $x,$  producing functions of $\mu$ alone.  This is as it should be, as the functions $g_{\o}(\mu) $ live in $x$-space but move with $\mu.$  
\bigskip

%Since $f^{|\o|}(g_{\o}(\mu))=c,$ where $|\o|=n,$ we agree to call the functions $g_{\o}(\mu)$ {\it level functions} of order $n,$  or $n^{th}$-order level functions corresponding to the word $\o.$
%\bigskip

Notice also that $f^{k}(g_{\o}(\mu))=g_{\s^k(\o)}(\mu),$ and in particular, $f^{|\o|}(g_{\o}(\mu))=c.$ 
\bigskip

%It is useful to define $g_c(\mu)=c$ for al $\mu.$

%For all $n\ge 0,$ let $G_n(\mu)$ denote the graph of $f_\mu^n.$
%
%By the chain rule, $\dfrac{d}{dx}f^n_\mu(x)=\prod_{k=0}^{n-1}\mu f'[f_\mu^k(x)].$
%Therefore,  if $x$ is an extreme point of $G_n(\mu)$ then there exists $k,\ 0\le k\le n$ such that $f^k(x)=c$.
%\bigskip
%
%{\bf Definition}\ \quad
%For all $\mu$ fixed, define $g_c(\mu)=c,$ and for $1\le k\le n-1,$ denote by the symbol $g_\o(\mu)$ any $k^{th}$ preimage of $c,$ $|\o|=k,$ specifically
%$$
%f^{-k}_\mu(c)=\{g_{\o_{1}\o_{2}\cdots \o_{k}}(\mu) \ \o=\o_1 \o_2 \dots \o_k\quad \o_i\in\{L,R\}\}.
%$$
%
%Fixed points of the  functions $g_{\o}$, will be central in the sequel. 

%Clearly, if $\u<c$ then $f^{-k}_\mu(c)$ is empty since $c$ has no preimages at all. So, for each sequence $P$ of $R's$ and $L's$ there are parameter values $\u$ for which $\xss{\o}{k}$ is undefined

%For example, for the map $\mu f(x)=4\mu x (1-x)$ we have
%$$
%(\mu f)^{-1}_L(x)=\frac{u-\sqrt{u^2-u x}}{2 u}
%$$
%\indent and
%$$
%(\mu f)^{-1}_R(x)=\frac{u+\sqrt{u^2-u x}}{2 u}\cdot
%$$
%Then, for example, with $P=RL$ and then $P=RLR$ 
%$$
%x^2_{RL}(\mu)=((\mu f)^{-1}_R) \circ (\mu f)^{-1}_L)(x)|_{x=.5}=\frac{u+\sqrt{u^2+\frac{1}{2} \left(\sqrt{u^2-0.5 u}-u\right)}}{2 u},
%$$
%$$
%x^3_{RLR}(\mu)=((\mu f)^{-1}_R) \circ (\mu f)^{-1}_L) \circ (\mu f)^{-1}_R)(x)|_{x=.5}=\frac{u+\sqrt{u^2+\frac{1}{2} \left(\sqrt{u^2+\frac{1}{2} \left(-u-\sqrt{u^2-0.5 u}\right)}-u\right)}}{2 u},
%$$
%which are purely functions of $\mu$ alone.
%\bigskip

For any $C^3$ function $\psi$, let $S(\psi)= \dfrac{\psi'''}{\psi'}-\dfrac{3}{2}\left(\dfrac{\psi''}{\psi'}\right)^2.$ 
A simple computation reveals that if $\phi$ is also a $C^3$ function and if $S(\phi)>0$, $S(\psi)>0$, then $S(\psi\circ \phi)>0.$

\section*{main section}

Here we prove that for each member of $\cal{C}$, as defined just below, the map $\mu\mapsto K(\mu f)$ is monotone, and receive as a corollary that the topological entropy $\mu\mapsto h(\mu f)$ is also monotone.
\bigskip

Consider the family $\cal{C}$ with the following properties:
\bigskip

1)\ For each $\mu$ there exists a unique fixed point for $\mu f$ in
$(0,1)$.  
\bigskip

2)\ For each fixed $\mu$ and for all $n\ge 1,$ $f^n_\mu$ has at most one attracting periodic orbit,
and $O(\u)$ is asymptotic to this attracting periodic orbit.
\bigskip

3)\  $S[(\mu f)_{\o}^{-1}]>0$  for all $\mu,$ where $\o\in\{L, R\}.$
\vskip 20pt

{\bf Remarks}
\bigskip

$\cal{C}\ne\emptyset$ since
\bigskip

(i)\ Concave maps, for example, have property 1.
\bigskip

(ii)\ It is known that if $S(f)<0$ for all $x$, then property 2 holds. [4]
\bigskip

(iii)\  One can check that $S[(\mu f)_{\o}^{-1}]>0$ for all $\mu$ ($\o\in\{L,R\})$ and that $S(f)<0$ for all $x$ when $f(x)=4x(1-x).$
\bigskip

%%%%%%%%%%%%%%%%%%%%%%%%%%%%%%%% definition of K was here %%%%%%%%%%%%%%%%%%%%%%%%

\bigskip
{\bf Lemma} Assume that there is a $\mu^*\in (c,1]$ with $g_{\o}(\mu^*)=\mu^*,$ where $\o=\o_1\o_2\cdots \o_{k-1}\in\cal{K}.$  Then
$f_{\u^*}^k(c)=c,$ that is, $\mu^*$ corresponds to a super stable $k$-cycle.
\smallskip

 In the proof the generic $\mu$ repalces the $\mu^*$ of the hypothesis for ease of reading.
\bigskip

%\begin{figure}[h] % float placement: (h)ere, page (t)op, page (b)ottom, other (p)age
%  \centering
%  % file name: C:/Documents and Settings/john/Desktop/TeXfiles/z.pdf
%  \includegraphics[width=2in,height=2in,keepaspectratio]{z}
%  \caption{Counter example to the converse of the lemma}
%  \label{}
%\end{figure}
%htbp bb=125 73 381 329,

{\bf Proof}\quad  $g_{\o}({\mu})=({\mu} f)_{\o_{1}}^{-1}\c\cdots\c ({\mu} f)_{\o_{k-1}}^{-1}(c)={\mu}\ \Rightarrow $
$$c=f^{k-1}_{\mu}\left[
g_\o({\mu})
\right]
=f^{k-1}_{\mu}\left[
({\mu} f)_{\o_{1}}^{-1}\c\cdots\c ({\mu} f)_{\o_{k-1}}^{-1}(c)
\right]
=f^{k-1}_{\mu}({\mu})=f^{k-1}_{\mu}\left[
{\mu} f(c)
\right]=f_{\mu}^k(c).$$
$\quad\square$
\bigskip

Note that by continuity there is an open set of parameter values containing $\mu^*$ for which the composition
$g_{\o}(\mu)=(\mu f)_{\o_1}^{-1}\circ \cdots\circ (\mu f)_{\o_{k-1}}^{-1}(c)$
is defined.
%On the other hand, 
%
%  Applying $(\u f)_{L}^{-1}$ or $(\u f)_{R}^{-1}$ along the word $\o$ to both sides of the equation $f_{\u}^k(c)=c$ gives
%$$x^{k-1}_{\o}({\mu})=(\u f)_{\\\o_{1}}^{-1}\c\cdots\c (\u f)_{\\\o_{k-1}}^{-1}(c)=\mu f(c)=\u.$$\quad$\square$
\bigskip

%{\bf Remarks}
%\bigskip
%
%(i)\ A counter example to the converse of the lemma is provided in figure 1, where f has a 6-cycle with kneading sequence $RlR^{\infty}.$  The kneading sequence has a two-fold splitting [9], or is not `tight' [10].
%\bigskip
%
%(ii)\ The trajectories of distinct preimages can never intersect.
\bigskip

\noindent {\bf Remarks:}\\ 
 Denote by $G_{n}(\mu)$ the graph of $f_\mu^n.$  It follows from the implicit function theorem that, for all $n\ge 1,\ 0\le k\le n-1,$ level functions of order $k$ exist so long as the intersection of $G_{n}(\mu)$ with the line $y=c$ exists.

But this intersection exists for all $\mu>\mu^*$, where $\mu^*$ is the parameter value with the property that, for $1\le k\le n,$ 
$g_{\o}(\mu^*)=\mu^*$,  $\o\in\K{k-1}$;  for then, $f^{k}_{\mu^*}(c)=c $ in $G_{n}(\mu)$ by the lemma.  Therefore, so long as $\mu^*$ is unique with the above property, we see that for all $\mu>\mu^*$, the intersection of $G_{n}(\mu)$ and the line $y=c$ persists, and so, the level functions $g_\o,\ \o\in\K{n-1}$ exist on a connected domain.

As mentioned above, a certain number of these $g_\o,\ \o\in\K{n-1}$ will have fixed points, and these will correspond to super stable points of period $n$ by the lemma.  
\bigskip

{\bf Theorem}\quad For each member of $\C$, the map $\mu\mapsto K(\mu f)$ is monotone.
\bigskip

{\bf Proof}\quad Since each member of $\C$ is of the form $\mu f,$ there exists a unique parameter value $\mu^*$, namely $\mu^*=c,$ such that for all $\mu>\mu^*,$ $\g{R}$ exists on the domain $(\u^*,1].$  $g_R$ is continuous on the connected set $(\u^*,1]$ and so the image of $g_R=$ im($g_R$) is connected.  But $S(g_R)>0$ by assumption, so that $\dfrac{d}{d\mu}g_R(\mu)$ cannot have a positive local maximum.  Therefore, there exists a unique $\mu$ such that $g_R(\mu)=\mu,$ or, by the lemma, a unique $\mu$ such that $f_\mu^2(c)=c,$ where $K(\mu f)=RC.$

By way of (strong) induction, assume that for some $n>1$ and for all $k, 1\le k\le n-1,$ and for all $\o\in\K{k}$ there is a unique $\mu_\o$ such that $g_{\o}(\mu_\o)=\mu_\o,$ or $f_{\mo}^k(c)=c.$ 
In particular, for every $\o\in\K{n-1},f_{\mo}^n(c)=c$ for a unique 
$\mo,$ where $K(\mu f)=\o_1\o_2\cdots\o_{n-1}C.$

Focusing on the domains of the functions $g_\o$, where $\o\in\cup_{k=1}^{n-1}\K{k},$ notice that by the induction hypothesis, there exists a unique $\mo$ such that $g_{\o}(\mu_\o)=\mu_\o$ which implies that for every $\epsilon >0,\ \mo+\epsilon\in$ dom$(g_{\o\tau}), \tau\in\{L,R\}.$   Therefore, dom$(g_{\o\tau})=(\mo,1]=\{\mu\le 1\ |\ g_\o(\mu)<\mu\}$ is connected, and as $g_{\o\tau}$ is a continuous function on this connected domain, im($g_{\o\tau})$ is connected.  

Now we claim that there exists $\delta>0,$ such that for every $\mu\in(\mo,\mo+\delta), g_{\o\tau}(\mu)>\mu.$  If we assume instead that there is no such $\delta,$  then there exist $\mu_1, \mu_2\in$ dom$(g_{\o\tau})$ with $g_{\o\tau}(\mu_1)<\mu_1$
and $g_{\o\tau}(\mu_2)\ge\mu_2,$ and by the intermediate value theorem, there exists $\mu',\ \mu_1<\mu'<\mu_2$ such that $g_{\o\tau}(\mu')=g_{\o}(\mu').$  But $f_{\mu'}^n[g_{\o\tau}(\mu')]=f_{\mu'}^n[g_{\o}(\mu')]$ implies that $c=\mu' f(c)=\mu'.$  Recall that $f(c)=1.$  This contradiction establishes the claim.

With this in mind, since dom$(g_{\o\tau})$ is connected and for each $\mu\in$ dom$(g_{\o\tau}),\ S(g_{\o\tau})>0,$ for each $\o\in\K{n}$ there exists a unique $\mo$ such that $g_{\o\tau}(\mo)=\mo,$ or, by the lemma, $f_{\mo}^{n+1}(c)=c,$ with $K(\mu f)=\o_1\o_2\cdots\o_{n}C.$

This concludes the argument and we have shown that for each member of $\C$, the map $\mu\mapsto K(\mu f)$ is monotone.  $\square$
\bigskip

{\bf Corollary}\quad For each member of $\C$, the map $\mu\mapsto h(\mu f)$ is monotone, where $h$ represents the topological entropy.
\bigskip

{\bf Proof}\quad
As the rate of orbit production for members of $\cal{C}$ can never decrease, the topological entropy of maps in the class $\cal{C}$ is evidently monotone with the parameter.  $\square$

.

% The proof is by strong induction.  First, notice that $\xss{R}{1}$ exists on the connected domain $[c,1]$.  Since we assume that $S(\xss{R}{1})>0$, $\du{\mu}\xss{R}{1}$ cannot have a positive local maximum.  Therefore,
%$$
%\e{!\us}\left[\xsss{R}{1}=\us\right]\ \Rightarrow\  \e!\us\left[f_{\us}^2(c)=c\right]\quad\text{for}\quad K(\us f)=RC\quad \text{by the lemma}.
%$$
%In other words, $\e!\mu_{\o}\ [f^1_{\mu_{\o}}(c)=c]$ with $\o=K(\mu_\o f) \raw\  \e!\mu_{\tau}\ [f^2_{\mu_{\tau}}(c)=c]$ with $\tau=K(\mu_{\tau} f).$  Here $\o=C$ and $\tau=RC$.
%\smallskip
%
%Assume that for $1\le k\le n,$ and for all $\o=K(\mu f)$ (with the length of $\o$ not exceeding $n$),  $\e{!\mu_\o}\ [x_{\o}^{k-1}(\mu_\o)=\mu_\o]$, that is $\e!\mu_\o\ [f_{\mu_\o}^k(c)=c]$ with $\o=K(\mu f).$
%\smallskip
%
%Since $\e{!\mu_\o}\ [x_{\o}^{k-1}(\mu_\o)=\mu_\o]$, that is, $\e!\mu_\o\ [f_{\mu_\o}^k(c)=c],$ dom$(\xss{\o}{k})$ is connected.
%% keeping in mind that, for all $P$ such that $P=K(\mu f),$
%%$$
%%\fa{\ep,\ 0<\ep\ll1},\ x_{Q}^{k}(\mu^*+\ep)<\mu^*+\ep<x_{\o}^{n}(\mu^*+\ep).
%%$$
%\smallskip
%
%If  $\xss{\o}{k}$ has a fixed point, that is, if $P=K(\mu f)$ for some $\mu$, then $f_{\mu}^{k+1}(c)=c$ when $\xss{\o}{k}=\mu.$
%\smallskip
%
%But $S[\xss{\o}{k}]>0\ \raw\ \e{!\mu}\ [x_{\o}^{k}(\mu)=\mu],$ that is, $\e!\mu\ [f_{\mu}^{k+1}(c)=c]$ with $\o=K(\mu f).$
%\smallskip
%
%In particular, $\e!\mu_\o\ [f_{\mu_\o}^n(c)=c]$ with $\o=K(\mu_{\o} f) \raw\  \e!\mu_{\tau}\ [f_{\mu_{\tau}}^{n+1}(c)=c]$ with $\tau=K(\mu_{\tau} f).$\quad$\square$
 \bigskip

{\bf Note:}\

One computes that $S[(\mu f)^{-1}_\o]>0,\ \o\in\{L,R\}$ when $\mu f(x)=4\mu x(1-x)$ and $\mu f(x)=\mu\sin(\pi x).$
\bigskip 

\end{document}